\documentclass[11pt, reqno]{amsart}
\usepackage{indentfirst, amssymb, amsmath, amsthm, mathrsfs, setspace, indentfirst, enumerate,  mathrsfs, amsmath, amsthm}
\usepackage[bookmarksnumbered, colorlinks, plainpages]{hyperref}
\usepackage{mathrsfs}
\newtheorem*{theo1A}{Theorem 1.A}
\newtheorem*{theo1B}{Theorem 1.B}

\newtheorem*{op1A}{Open Problem 1.A}

\newtheorem{ques}{Question}[section]

\newtheorem*{cor A}{Corollary A}
\newtheorem*{cor B}{Corollary B}

\newtheorem{theo}{Theorem}[section]
\newtheorem{lem}{Lemma}[section]
\newtheorem{cor}{Corollary}[section]

\newtheorem{exm}{Example}[section]

\newcommand{\ol}{\overline}

\newcommand{\be}{\begin{equation}}
\newcommand{\ee}{\end{equation}}
\newcommand{\beas}{\begin{eqnarray*}}
\newcommand{\eeas}{\end{eqnarray*}}
\newcommand{\bea}{\begin{eqnarray}}
\newcommand{\eea}{\end{eqnarray}}

\numberwithin{equation}{section}
\begin{document}

\title[S\MakeLowercase{olutions of Certain} F\MakeLowercase{ermat-Type Partial} ......]{\LARGE S\Large\MakeLowercase{olutions of Certain} F\MakeLowercase{ermat-Type Partial Differential-Difference equations}}

\date{}
\author[S. M\MakeLowercase{ajumder} \MakeLowercase{and} D. P\MakeLowercase{ramanik}]{S\MakeLowercase{ujoy} M\MakeLowercase{ajumder}$^{*}$ \MakeLowercase{and} D\MakeLowercase{ebabrata} P\MakeLowercase{ramanik}}
\address{Department of Mathematics, Raiganj University, Raiganj, West Bengal-733134, India.}
\email{sm05math@gmail.com, sjm@raiganjuniversity.ac.in}
\address{Department of Mathematics, Raiganj University, Raiganj, West Bengal-733134, India.}
\email{debumath07@gmail.com}

\renewcommand{\thefootnote}{}
\footnote{*\emph{Corresponding Author}: Sujoy Majumder.}

\renewcommand{\thefootnote}{\arabic{footnote}}
\setcounter{footnote}{0}

\begin{abstract} The purpose of this paper is to investigate the non-constant entire as well as meromorphic solutions of the Fermat-type partial differential-difference equation:
\[\left(\sum_{j=1}^m\frac{\partial f(z_1, z_2, \ldots, z_m)}{\partial z_j}\right)^{m_1} + f^{m_2}(z_1 + c_1, z_2 + c_2, \ldots, z_m + c_m ) = 1,\]
where $m_1$ and $m_2$ are positive integers such that $m_1+m_2>2$ and $(c_1, c_2, \ldots, c_m)\in \mathbb{C}^m$. The results of our paper generalize the result of Xu and Wang \cite {XW1} from $\mathbb{C}^2$ to $\mathbb{C}^m$. Also in the paper we give positive answer of the open problem addressed by Xu and Wang \cite {XW1}. Moreover plenty of examples are provided to illustrate our findings.
\end{abstract}

\thanks{Typeset by \AmS -\LaTeX}
\maketitle

\section{{\bf Introduction and main results}}
We define $\mathbb{Z}_+=\mathbb{Z}[0,+\infty)=\{n\in \mathbb{Z}: 0\leq n<+\infty\}$ and $\mathbb{Z}^+=\mathbb{Z}(0,+\infty)=\{n\in \mathbb{Z}: 0<n<+\infty\}$. As on $\mathbb{C}^m$, the exterior derivative $d$ splits $d= \partial+ \bar{\partial}$ and twists to $d^c= \frac{i}{4\pi}\left(\bar{\partial}- \partial\right)$ (see \cite{HLY1,WS1}). Clearly $dd^{c}= \frac{i}{2\pi}\partial\bar{\partial}$.
An exhaustion $\tau_m$ of $\mathbb{C}^m$ is defined by $\tau_m(z)=||z||^2$. The standard Kaehler metric on $\mathbb{C}^m$ is given by $\upsilon_m=dd^c\tau_m>0$. On $\mathbb{C}^m\backslash \{0\}$, we define $\omega_m=dd^c\log \tau_m\geq 0$ and $\sigma_m=d^c\log \tau_m \wedge \omega_m^{m-1}$. For any $S\subseteq \mathbb{C}^m$, let $S[r]$, $S(r)$ and $S\langle r\rangle$ be the intersection of $S$ with
respectively the closed ball, the open ball, the sphere of radius $r>0$ centered at $0\in \mathbb{C}^m$.

\smallskip
We define 
\[\displaystyle \partial_{z_i}(f(z))=\frac{\partial f(z)}{\partial z_i},\ldots,\partial^{l_i}_{z_i}(f(z))=\frac{\partial^{l_i} f(z)}{\partial z_i^{l_i}},\]
where $l_i\in \mathbb{Z}^m_+$ and $i=1,2,\ldots,m$ and 
\[\partial^{I}(f(z))=\frac{\partial^{|I|}f(z)}{\partial z_1^{i_1}\cdots \partial z_m^{i_m}},\]
where $|I|=\sum_{j=1}^m i_j$. Therefore the zero multiplicity $\mu^0_f(a)$ of $f$ at a point $a\in \mathbb{C}^m$ is defined
to be the order of vanishing of $f$ at $a$. If $a=(a_1,\ldots,a_m)$, then $\partial^I(f(a))=0$, where $|I|\leq \mu^0_f(a)-1$.
In other words, we can write $f(z)=\sum_{i=0}^{\infty}P_i(z-a)$, where the term $P_i(z-a)$ is either identically zero or a homogeneous polynomial of degree $i$. Certainly $\mu^0_f(a)=\min\{i:P_i(z-a)\not\equiv 0\}$.

\smallskip
Let $f$ be a meromorphic function on $G$. Then there exist holomorphic functions $g$ and $h$ such that $hf=g$ on $G$ and $\dim_z h^{-1}(\{0\})\cap g^{-1}(\{0\})\leq m-2$. Therefore the $c$-multiplicity of $f$ is just $\mu^c_f=\mu^0_{g-ch}$ if $c\in\mathbb{C}$ and $\mu^c_f=\mu^0_h$ if $c=\infty$. The function $\mu^c_f: \mathbb{C}^m\to \mathbb{Z}$ is nonnegative and is called the $c$-divisor of $f$. If $f\not\equiv 0$ on each component of $G$, then $\nu=\mu_f=\mu^0_f-\mu^{\infty}_f$ is called the divisor of $f$. We define 
$\text{supp}\; \nu=\text{supp}\;\mu_f=\ol{\{z\in G: \nu(z)\neq 0\}}$.

\smallskip
Let $f$, $g$ and $a$ be meromorphic functions in $\mathbb{C}^m$. Then one can find three pairs
of entire functions $f_1$ and $f_2$, $g_1$ and $g_2$, and $a_1$ and $a_2$, in which each pair is coprime
at each point in $\mathbb{C}^m$ such that $f = f_2/f_1$, $g=g_2/g_1$ and $a = a_2/a_1$.
Now $f$ and $g$ share $a$ CM if $\mu_{a_1f_2-a_2f_1}^0=\mu_{a_1g_2-a_2g_1}^0\;(a\not\equiv \infty)$ and $\mu_{f_1}^0=\mu_{g_1}^0\;\;(a=\infty)$. 
For $t>0$, the counting function $n_{\nu}$ is defined by
\beas n_{\nu}(t)=t^{-2(m-1)}\int_{A[t]}\nu \upsilon_m^{m-1},\eeas
where $A=\text{supp}\;\nu$. 
The valence function of $\nu$ is defined by 
\[N_{\nu}(r)=N_{\nu}(r,r_0)=\int_{r_0}^r n_{\nu}(t)\frac{dt}{t}\;\;(r\geq r_0).\]

Also we write $N_{\mu_f^a}(r)=N(r,a;f)$ if $a\in\mathbb{C}$ and $N_{\mu_f^a}(r)=N(r,f)$ if $a=\infty$.
For $k\in\mathbb{N}$, define the truncated multiplicity functions on $\mathbb{C}^m$ by $\mu_{f,k}^a(z)=\min\{\mu_f^a(z),k\}$, $\mu_{f(k}^a(z)=\mu_f^a(z)$ if $\mu_f^a(z)\geq  k$ and $\mu_{f(k}^a(z)=0$ if $\mu_f^a(z)<k$. Also we define
the truncated valence functions $N_{\nu}(t)=\ol N(t,a;f)$ if $\nu=\mu_{f,1}^a$ and $N_{\nu}(t)=N_{(k}(t,a;f)$, if $\nu=\mu_{f(k}^a$.

An algebraic subset $X$ of $\mathbb{C}^m$ is defined as a subset
\[X=\left\lbrace z\in\mathbb{C}^m: P_j(z)=0,\;1\leq j\leq l\right\rbrace\]
with finitely many polynomials $P_1(z),\ldots, P_l(z)$.
A divisor $\nu$ on $\mathbb{C}^m$ is said to be algebraic if $\nu$ is the zero divisor of a polynomial. In this case the counting
function $n_{\nu}$ is bounded (see \cite{GK1,WS1}).

\medskip
With the help of the positive logarithm function, we define the proximity function of $f$ by
\[m(r, f)=\mathbb{C}^m\langle r; \log^+ | f | \rangle=\int_{\mathbb{C}^m\langle r\rangle} \log^+ |f|\;\sigma_m.\]

The characteristic function of $f$ is defined by $T(r, f)=m(r,f)+N(r,f)$. We define $m(r,a;f)=m(r,f)$ if $a=\infty$ and $m(r,a;f)=m(r,1/(f-a))$ if $a$ is finite complex number. Now if $a\in\mathbb{C}$, then the first main theorem of Nevanlinna theory states that $m(r,a;f)+N(r,a;f)=T(r,f)+O(1)$, where $O(1)$ denotes a bounded function when $r$ is sufficiently large. 
We define the order of $f$ by
\[\rho(f):=\limsup _{r \rightarrow \infty} \frac{\log T(r, f)}{\log r} .\]

We define the linear measure $m(E):=\int_E dt$, the logarithmic measure $l(E):=\int_{E\cap [1,\infty)} \frac{d t}{t}$ and the upper density measure
\[\ol {\text{dens}}\;E=\lim\limits_{r\rightarrow \infty}\frac{1}{r}\int_{E\cap [1,r]} dt\]
for a set $E\subset [0,\infty)$. Moreover, if $l(E)<+\infty$, resp., $m(E)<+\infty$, then $E$ is of zero upper density.
Let $S(f)=\{g:\mathbb{C}^m\to\mathbb{P}^1\;\text{meromorphic}: \parallel T(r,g)=o(T(r,f))\}$, where $\parallel$ indicates that the equality holds only outside a set $E$ with zero upper density measure and the element in $S(f)$ is called the small function of $f$.

\subsection{{\bf Fermat type functional equation in $\mathbb{C}$}} The following functional equation
\bea\label{ne1} f^n(z)+g^n(z)=1,\eea
where $n$ is a positive integer can be regarded as the Fermat diophantine equations $x^n+y^n=1$ over function fields. 
In 1927, Montel \cite{M1} proved that the equation (\ref{ne1}) has no transcendental entire solutions for $n \geq 3$.
Also in Theorem 4 \cite{FG2}, we see that the equation (\ref{ne1}) actually has no non-constant entire solutions when $n>2$.
Gross \cite[Theorem 3]{FG1} proved that equation (\ref{ne1}) has no non-constant meromorphic solutions when $n>3$.  
For $n=2$, Gross \cite[Theorem 4]{FG2} found that the equation (\ref{ne1}) has entire solutions of the from $f(z)=\sin (h(z))$ and $g(z)=\cos (h(z))$, where $h(z)$ is an entire function.
For $n=2$ and $n=3$, Gross \cite[Theorem 1]{FG1} and Baker \cite[Theorem 1]{Baker} respectively obtained the solutions of Eq. (\ref{ne1}) and it specific forms.
In 1970, Yang \cite{Y1} investigated the following Fermat-type equation
\bea\label{ne2} f^{m_1}(z)+g^{m_2}(z)=1\eea
and obtained that equation (\ref{ne2}) has no non-constant entire solutions when $\frac{1}{m_1}+\frac{1}{m_2}<1$ (see the proof of Theorem 1). Therefore it is clear that the equation (\ref{ne2}) has no non-constant entire solutions when $m_1>2$ and $m_2 > 2$. However, for the case when $m_1=m_2=2$ and $g(z)$ has a specific relationship with $f(z)$, many authors investigated the existence of solutions of the equation (\ref{ne2}). As a result, successively several research papers were published (see \cite{Liu1}{-\cite{Lu2}, \cite{TL1}, \cite{YL1}) and the research in this area is still ongoing.

\subsection{{\bf Fermat type functional equation in $\mathbb{C}^m$}}
In 2008, Li \cite[Theorem 1]{L1} proved that meromorphic solutions of the Fermat-type functional equation $f^2+g^2=1$ in $\mathbb{C}^2$ must be constant if and only if $\frac{\partial f(z_1,z_2)}{\partial z_2}$ and $\frac{\partial g(z_1,z_2)}{\partial z_1}$ share $0$ CM.  Moreover if $f(z_1,z_2)=\frac{\partial u(z_1,z_2)}{\partial z_1}$ and $g(z_1,z_2)=\frac{\partial u(z_1,z_2)}{\partial z_2}$, then any entire solutions of the partial differential equations:
\[\left(\frac{\partial u(z_1,z_2)}{\partial z_1}\right)^2+\left(\frac{\partial u(z_1,z_2)}{\partial z_2}\right)^2=1\]
in $\mathbb{C}^2$ are necessarily linear (see \cite{DK1}). From the fact that the surface $z_1^n+z_2^n= 1$ in $\mathbb{C}^2$ is a Kobayashi hyperbolic manifold and thus there are no non-constant holomorphic curves $(f_1, f_2)$ from the complex plane $\mathbb{C}$ (and thus $\mathbb{C}^m$) to the surface (see \cite[page 360]{BS1}), one can easily prove that any
entire solutions of the partial differential equations
\[\left(\frac{\partial u(z_1,z_2)}{\partial z_1}\right)^n+\left(\frac{\partial u(z_1,z_2)}{\partial z_2}\right)^n=1,\]
where $n\geq 3$ must also be linear. The solutions of Fermat-type partial differential equations
was originally investigated by Li \cite{L2} and Saleeby \cite{GS1}.

\medskip
Recently, Fermat-type difference and Fermat-type partial differential-difference equations are the current interests among researchers by utilizing Nevanlinna value distribution theory in several complex variables (see \cite{AA1}-\cite{AM1}, \cite{MSP}, \cite{MA}, \cite{GS2}, \cite{XW1}, \cite{XH}).

In 2020, Xu and Wang \cite {XW1} considered the following partial differential-difference equation:
\bea\label{f01} \left(\frac{\partial f(z_1, z_2)}{\partial z_1}+\frac{\partial f(z_1, z_2)}{\partial z_2}\right)^{m_1}+f^{m_2}(z_1+c_1, z_2+c_2)=1,\eea
to discuss the existence of finite order transcendental entire solutions, where $m_1$ and $m_2$ are positive integers such that $m_1+m_2>2$.  Using the difference Nevanlinna theory for meromorphic functions (see \cite{TBC1, CK1, CX1, K1}) in the case of higher dimension,  Xu and Wang \cite {XW1} obtained the following result.

\begin{theo1A}\cite[Theorem 1.6]{XW1}
Let $c=(c_1, c_2)\in\mathbb{C}^2$ and $m_1, m_2$ be two distinct positive integers such that $m_1+m_2>2$. If the equation (\ref{f01}) satisfies one of the conditions: $(i)$ $m_2>m_1$ $(ii)$ $m_1>m_2\geq 2$,
then equation (\ref{f01}) does not have any finite order transcendental entire solutions.
\end{theo1A}

When $m_1=2$ and $m_2=1$,  Xu and Wang \cite{XW1} obtained the following result regarding the finite order transcendental entire solutions of the equation (\ref{f01}).

\begin{theo1B}\cite[see pp.10]{XW1}
Let $c=(c_1, c_2)\in\mathbb{C}^2$. If $f(z_1,z_2)$ is a finite order transcendental entire solution of the following equation:
\beas\label{f03} \left(\frac{\partial f(z_1, z_2)}{\partial z_1}+\frac{\partial f(z_1, z_2)}{\partial z_2}\right)^{2}+f(z_1+c_1, z_2+c_2)=1,\eeas
then $f(z_1,z_2)$ takes the from
\beas &&f(z_1,z_2)\\&=&1-\frac{1}{4}c_1^2-\frac{1}{4}z_1^2+z_1\left\lbrack G_2(z_2-z_1)+a_3(z_2-z_1)\right\rbrack -c_1G(z_2-z_1)-a_3c_1[z_2-z_1-(c_2-c_1)]-\\&&
\left(G_2(z_2-z_1)+a_3 \left\lbrack z_2-z_1 -(c_2-c_1)\right\rbrack \right)^2 ,\eeas
where $G_2(s)$ is a finite order transcendental entire period function with period $c_2-c_1$ and $a_3=\frac{c_1}{2(c_2-c_1)}$.
\end{theo1B}

In the same paper, Xu and Wang \cite{XW1} imposed the following open problem for further investigation:
\begin{op1A}\label{open} Whether there exists the finite order transcendental entire solutions for the equation (\ref{f01}) in the case $m_1>2$ and $m_2=1$ or not?
\end{op1A}

To our knowledge, Open Problem 1.A remain unresolved up to the present. Consequently our first goal is to give affirmative answer of 
the open problem 1.A. In the paper, we consider the following Fermat-type partial differential-difference equation
\bea\label{pds1}F^{m_1}(z) + f^{m_2}(z+c)=1,\eea
where 
\bea\label{lm0} F(z)=\sideset{}{_{j=1}^m}{\sum}\frac{\partial f(z_1, z_2, \ldots, z_m)}{\partial z_j}.\eea
Moreover, for $m_1=m_2=1$, it is easy to verify that $f(z)=e^{z_1+z_2+\cdots+z_m}+1$ is a solution of 
$F(z)+f(z+c)=1$, where $e^{c_1+c_2+\cdots+c_m}=-m$. Therefore we consider the Eq. (\ref{pds1}) for the existence of solutions for the case when $m_1+m_2>2$. 

Our next goal is to investigate all possible finite-order transcendental entire as well as meromorphic solutions of the Eq. (\ref{pds1}) for all different combinations of $m_1$ and $m_2$. Specifically, we focus on understanding the conditions under which such solutions exist and determining their explicit forms.  Furthermore, we extend these findings to higher-dimensional settings $\mathbb{C}^m$ and explore how the solutions behave when multiple variables are involved.
Now we state our main results.

\begin{theo}\label{t2.1} Let $c\in\mathbb{C}^m$ and let $m_1$ and $m_2$ be two positive integers such that $m_1+m_2>2$. Then any finite order entire solutions of the equation (\ref{pds1}) are characterized as follows:
\begin{enumerate}
\item[(i)] if $m_2\geq m_1$, then $m_1=m_2=2$ and $f(z)=\sin (d_1z_1+d_2z_2+\ldots+d_mz_m+P(z))$,
where $d_1,d_2,\ldots,d_m\in \mathbb{C}$ and $P(z)$ is a polynomial in $\mathbb{C}^m$ such that $P(z+c)\equiv P(z)$ and
\beas e^{\iota (d_1c_1+d_2c_2+\ldots+d_mc_m)}\sideset{}{_{j=1}^m}{\sum} \left(\frac{\partial P(z)}{\partial z_j}+d_j\right)=1,\eeas
\item[(ii)] if $m_1>m_2$, then $m_1=2$, $m_2=1$ and 
\[f(z)=1-\left(-(z_1+z_2+\ldots+z_m)/2m+g(z)\right)^2,\]
where $g(z)=g(z_2-z_1,z_3-z_2,\ldots,z_m-z_{m-1})$ is an arbitrary finite order transcendental entire function such that 
$g(z+c)-g(z)=(c_1+c_2+\ldots+c_m)/2m$.
\end{enumerate}
\end{theo}
 
Regarding the Open Problem 1.A, we get the following corollary from Theorem \ref{t2.1}.
\begin{cor} For $m_2=1$ and $m_1>2$, the Fermat-type partial differential-difference equation (\ref{pds1}) does not have any finite order transcendental entire solutions.
\end{cor}

If $c=(0,0,\ldots,0)$, then Theorem \ref{t2.1} implies the following corollary.
\begin{cor} Let $m_1$ and $m_2$ be two positive integers such that $m_1+m_2>2$. Then any non-constant entire solutions of the following equation 
\bea\label{ll1} F^{m_1}(z)+f^{m_2}(z)=1\eea
are characterized as follows:
\begin{enumerate}
\item[(i)] if $m_2\geq m_1$, then $m_1=m_2=2$ and $f(z)=\sin (d_1z_1+d_2z_2+\ldots+d_mz_m+P(z))$,
where $d_1,d_2,\ldots,d_m\in \mathbb{C}$ and $P(z)$ is a polynomial in $\mathbb{C}^m$ such
$\sum_{j=1}^m \left(\frac{\partial P(z)}{\partial z_j}+d_j\right)=1$,
\item[(ii)] if $m_1>m_2$, then $m_1=2$, $m_2=1$ and 
\[f(z)=1-\left(-(z_1+z_2+\ldots+z_m)/2m+g(z)\right)^2,\]
where $g(z)=g(z_2-z_1,z_3-z_2,\ldots,z_m-z_{m-1})$ is an arbitrary entire function.
\end{enumerate}
\end{cor}

Following examples are provided to illustrate our findings.
\begin{exm} Let $f(z_1, z_2, z_3) = \sin (d_1z_1+d_2z_2+d_3z_3+P(z))$, where $d_1=d_2=d_3=\frac{1}{3}$, $P(z)=z_1-2z_2+z_3$ and $c=(\pi,0,-\pi)$.
Note that $f(z_1, z_2, z_3)$ is a finite order transcendental entire function and $e^{\iota (d_1c_1+d_2c_2+d_3c_3)}\sum_{j=1}^m \left(\frac{\partial P(z)}{\partial z_j}+d_j\right)=e^0=1$. It is easy to verify that $f(z_1, z_2, z_3)$ is a solution of the following Fermat type differential-difference equation
\[\left(\sideset{}{_{j=1}^3}{\sum}\frac{\partial f\left(z_1, z_2, z_3\right)}{\partial z_j}\right)^2+f^2(z_1+\pi, z_2, z_3-\pi)=1.\]
\end{exm}

\begin{exm} Let $f(z_1, z_2, z_3) = 1 - \left(-\frac{1}{6}( z_1+z_2+z_3) +g(z)\right)^2,$
where $g(z)=e^{z_1-2z_2 + z_3}$ and $c=(1, 0, -1)$. Note that $g(z)$ is a transcendental entire function of finite order such that $g(z+c)-g(z)=0=\frac{1}{6}(c_1+c_2+c_3)$. It is easy to verify that $f(z_1, z_2, z_3)$ is a finite order transcendental entire solution of the following Fermat type differential-difference equation
\[\left(\sideset{}{_{j=1}^3}{\sum}\frac{\partial f\left(z_1, z_2, z_3\right)}{\partial z_j}\right)^2+f\left(z_1+1, z_2, z_3-1\right)=1.\]
\end{exm}

For the existence of finite order non-constant meromorphic solutions of Fermat-type partial differential-difference equation (\ref{pds1}), we get the following results.

\begin{theo}\label{t2.2} Let $c\in\mathbb{C}^m$ and let $m_1$ and $m_2$ be two positive integers such that $m_1+m_2>2$. If one of the following conditions:
\begin{enumerate}
\item[(i)] $m_2>m_1\geq 3$;
\item[(ii)] $m_2>2m_1$ for $m_1\in\{1,2\}$,
\end{enumerate}
then the equation (\ref{pds1}) does not have any finite order meromorphic solution.
\end{theo}

\begin{theo}\label{t2.3} Let $c\in\mathbb{C}^m$ and let $m_1$ and $m_2$ be two positive integers such that $m_1+m_2>2$ and $m_1\geq m_2$. Then any finite order non-constant meromorphic solution $f(z)$ of the equation (\ref{pds1}) such that $\frac{\partial f(z)}{\partial z_1}$ and $\sum_{j=1}^m\frac{\partial f(z)}{\partial z_j}$ share $\infty$ CM are characterized as follows:
\begin{enumerate}
\item[(i)] if $m_1=m_2$, then $m_1=m_2=2$ and 
\[f(z)=\frac{h(z-c)-1/h(z-c)}{2\iota},\]
where $h(z)$ is a non-constant finite order meromorphic function in $\mathbb{C}^m$ such that
\[\iota\left(h(z+c)+1/h(z+c)\right)=\sideset{}{_{j=1}^m}{\sum}\frac{\partial h(z)}{\partial z_j}\left(1+1/h^2(z)\right)\]
and
\[\parallel\;\ol N(r,h(z-c))+\ol N\left(r,0;h(z-c)\right)=o(T(r,h)).\]

In particular if $\parallel\;N_{(2}(r,f)=o(T(r,f))$, then 
\[f(z)=\sin (d_1z_1+d_2z_2+\ldots+d_mz_m+P(z)),\]
where $d_1,d_2,\ldots,d_m\in \mathbb{C}$ and $P(z)$ is a polynomial in $\mathbb{C}^m$ such that $P(z+c)\equiv P(z)$ and $e^{-\iota (d_1c_1+d_2c_2+\ldots+d_mc_m)}\sum_{j=1}^m\left(\frac{\partial P(z)}{\partial z_j}+d_j\right)=\pm 1$.
Moreover if $c=0$, then 
\[f(z)=\sin \left((z_1+z_2+\ldots +z_m)/m+P(z)\right),\]
where $P(z)=P(z_3-z_2,\ldots,z_m-z_{m-1})$ is a polynomial such that $\sum_{j=1}^m\frac{\partial P(z)}{\partial z_j}=0$;
\item[(ii)] if $m_1>m_2$, then $m_1=2$, $m_2=1$ and 
\[f(z)=1-\left(-(z_1+z_2+\ldots+z_m)/2m+g(z)\right)^2,\]
where $g(z)=g(z_2-z_1,z_3-z_2,\ldots,z_m-z_{m-1})$ is an arbitrary finite order transcendental entire function such that 
$g(z+c)-g(z)=(c_1+c_2+\ldots+c_m)/2m$.
\end{enumerate}
\end{theo}

If $c=(0,0,\ldots,0)$, then Theorem \ref{t2.3} implies the following corollary.
\begin{cor} Let $m_1$ and $m_2$ be two positive integers such that $m_1+m_2>2$ and $m_1\geq m_2$. Then any non-constant meromorphic solution $f(z)$ of the equation (\ref{ll1}) such that $\frac{\partial f(z)}{\partial z_1}$ and $\sum_{j=1}^m\frac{\partial f(z)}{\partial z_j}$ share $\infty$ CM are characterized as follows:
\begin{enumerate}
\item[(i)] if $m_1=m_2$, then $m_1=m_2=2$ and $f(z)=\sin (d_1z_1+d_2z_2+\ldots+d_mz_m+P(z))$,
where $P(z)=P(z_3-z_2,\ldots,z_m-z_{m-1})$ is a polynomial such that $\sum_{j=1}^m\left(\frac{\partial P(z)}{\partial z_j}+d_j\right)= 1$.
\item[(ii)] if $m_1>m_2$, then $m_1=2$, $m_2=1$ and 
\[f(z)=1-\left(-(z_1+z_2+\ldots+z_m)/2m+g(z)\right)^2,\]
where $g(z)=g(z_2-z_1,z_3-z_2,\ldots,z_m-z_{m-1})$ is an arbitrary entire function.
\end{enumerate}
\end{cor}

If the set of poles of the solution $f(z)$ of the equation (\ref{pds1}) is algebraic, then we get the following result.

\begin{theo}\label{t2.4} Let $c\in\mathbb{C}^m$ and let $m_1$ and $m_2$ be two positive integers such that $m_1+m_2>2$. Then any finite order non-constant meromorphic solution $f(z)$ of the equation (\ref{pds1}) such that the set of poles of $f(z)$ is algebraic are characterized as follows:
\begin{enumerate}
\item[(i)] if $m_2\geq m_1$, then $m_1=m_2=2$ and $f(z)=\sin (d_1z_1+d_2z_2+\ldots+d_mz_m+P(z))$,
where $d_1,d_2,\ldots,d_m\in \mathbb{C}$ and $P(z)$ is a polynomial in $\mathbb{C}^m$ such that $P(z+c)\equiv P(z)$ and $e^{-\iota (d_1c_1+d_2c_2+\ldots+d_mc_m)}\sum_{j=1}^m\left(\frac{\partial P(z)}{\partial z_j}+d_j\right)=\pm 1$;

\item[(ii)] if $m_1>m_2$, then $m_1=2$, $m_2=1$ and 
\[f(z)=1-\left(-(z_1+z_2+\ldots+z_m)/2m+g(z)\right)^2,\]
where $g(z)=g(z_2-z_1,z_3-z_2,\ldots,z_m-z_{m-1})$ is an arbitrary finite order transcendental entire function such that 
$g(z+c)-g(z)=\frac{1}{2m}(c_1+c_2+\ldots+c_m)$.
\end{enumerate}
\end{theo}

Now based on Theorems \ref{t2.2} and \ref{t2.3}, we ask the following questions: 
\begin{ques} What can be said for non-constant meromorphic solutions of the Fermat-type partial differential-difference equation (\ref{pds1}) when $m_2>m_1\geq 1$? 
\end{ques}

\begin{ques} Is it possible to establish Theorem \ref{t2.3} without the hypothesis ``$\frac{\partial f(z)}{\partial z_1}$ and $\sum_{j=1}^m\frac{\partial f(z)}{\partial z_j}$ share $\infty$ CM''?
\end{ques}

\section{\bf{Auxiliary lemmas}}
The following lemmas are relevant to this paper and will be used in the subsequent analysis.
First we recall the lemma of logarithmic derivative:

\begin{lem}\label{L.1} \cite[Lemma 1.37]{HLY1} Let $f$ be a non-constant meromorphic function in $\mathbb{C}^m$ and $I=(\alpha_1,\ldots,\alpha_m)\in \mathbb{Z}^m_+$ be a multi-index. Then for any $\varepsilon>0$, we have
\[m\left(r,\partial^I(f)/f\right)\leq |I|\log^+T(r,f)+|I|(1+\varepsilon)\log^+\log T(r,f)+O(1)\]
for all $r\not\in E$, where $l(E)<+\infty$.
\end{lem}

The following result is known as second main theorem:
\begin{lem}\label{L.2} \cite[Lemma 1.2]{HY1} Let $f$ be a non-constant meromorphic function in $\mathbb{C}^m$ and let $a_1,a_2,\ldots,a_q$ be different points in $\mathbb{C}\cup \{\infty\}$. Then
\beas \parallel (q-2)T(r,f)\leq \sideset{}{_{j=1}^{q}}{\sum} \ol N(r,a_j;f)+O(\log (rT(r,f))).\eeas
\end{lem}

\begin{lem}\label{L.3} \cite[Theorem 1.26]{HLY1} Let $f$ be non-constant meromorphic function in $\mathbb{C}^m$. Assume that 
$R(z, w)=\frac{A(z, w)}{B(z, w)}$. Then
\beas T\left(r, R_f\right)=\max \{p, q\} T(r, f)+O\Big(\sideset{}{_{j=0}^p}{\sum} T(r, a_j)+\sideset{}{_{j=0}^q}{\sum}T(r, b_j)\Big),\eeas
where $R_f(z)=R(z, f(z))$ and two coprime polynomials $A(z, w)$ and $B(z,w)$ are given
respectively as follows: $A(z,w)=\sum_{j=0}^p a_j(z)w^j$ and $B(z,w)=\sum_{j=0}^q b_j(z)w^j$.
\end{lem}

\begin{lem}\label{L.6}\cite[Lemma 3.2]{HLY1} Let $f_j\not\equiv 0\;(j=1,2,\ldots,n;n\geq 3)$ be meromorphic functions in $\mathbb{C}^m$ such that $f_1,\ldots,f_{n-1}$ are non-constants and $f_1+\cdots+f_n=1$. If
\beas \parallel \;\;\sideset{}{_{j=1}^n}{\sum}\Big\lbrace N_{n-1}(r,0;f_j)+(n-1)\ol{N}(r,f_j)\Big\rbrace<\lambda T(r,f_j)+O(\log^+(T(r,f_j))\eeas
holds for $j=1,2,\ldots,n-1$, where $\lambda <1$, then $f_n\equiv 1$.
\end{lem}

\begin{lem}\label{L.8}\cite[Theorem 2.1]{CX1} Let $f$ be a non-constant meromorphic function in $\mathbb{C}^m$ and let $c\in \mathbb{C}^m$. If 
\[\lim\limits_{r\rightarrow \infty}\sup \frac{\log T(r,f)}{r}=0,\]
then
\beas \parallel\;m\left(r,f(z+c)/f(z)\right)+m\left(r,f(z)/f(z+c)\right)=o(T(r,f)).\eeas
\end{lem}

\begin{lem}\label{L.7}\cite[Theorem 2.2]{CX1} Let $f$ be a non-constant meromorphic function on $\mathbb{C}^m$ with 
\[\lim\limits_{r\rightarrow \infty}\sup \frac{\log T(r,f)}{r}=0,\]
then $\parallel\;T(r,f(z+c))=T(r,f)+o(T(r,f))$
holds for any constant $c\in\mathbb{C}^m$.
\end{lem}

Let $f$ be a non-constant meromorphic function on $\mathbb{C}^m$. Define complex differential-difference polynomials as follows:
\bea\label{cl1} P(f(z))=\sum\limits_{\mathbf{p} \in I} a_{\mathbf{p}}(z) f^{p_0}(z)\big(\partial^{\mathbf{i}_1}(f(z))\big)^{p_1} \cdots\big(\partial^{\mathbf{i}_l}(f(z))\big)^{p_l}f^{p_{l+1}}(z + \hat q_{l+1}) \cdots f^{p_{l+s}}(z + \hat q_{l+s}),\eea
$\mathbf{p}=\left(p_0, \ldots, p_{l+s}\right) \in \mathbb{Z}_{+}^{l+s+1}$,
\bea\label{cl2} Q(f(z))=\sum\limits_{\mathbf{q} \in J} c_{\mathbf{q}}(z) f^{q_0}(z)\big(\partial^{\mathbf{j}_1}(f(z))\big)^{q_1} \cdots\big(\partial^{\mathbf{j}_l}(f(z))\big)^{q_l}f^{q_{l+1}}(z+\tilde q_{l+1}) \cdots f^{q_{l+t}}(z+\tilde q_{l+t}),\eea
$\mathbf{q}=\left(q_0, \ldots, q_{l+t}\right) \in \mathbb{Z}_{+}^{l+t+1}$ and 
\bea\label{cl3} B(f(z))=\sideset{}{_{k=0}^{n}}{\sum} b_k(z) f^k(z),\eea
where $I, J$ are finite sets of distinct elements, $a_{\mathbf{p}}, b_k, c_{\mathbf{q}}\in S(f)$ such that $b_n\not\equiv 0$ and $\hat q_i, \tilde q_j\in\mathbb{C}^m$. 

For the case, when $P(f(z))$ and $Q(f(z))$ are complex differential polynomials, Hu and Yang \cite[Lemma 2.1]{ps1} generalised Clunie-lemma to high dimension. In 2020, Cao and Xu \cite[Theorem 3.6]{CX1} improved and extended Laine-Yang's difference analogue of Clunie theorem in one variable \cite[Theorem 2.3]{LY1} to high dimension by using Lemma \ref{L.8}. Now by using Lemmas \ref{L.1} and \ref{L.8} and proceeding in the same way as done in the proofs of Lemma 2.1 \cite{ps1} and Theorem 2.3 \cite{LY1} we get the following lemma.

\begin{lem}\label{L.4} Let $f(z)$ be a non-constant meromorphic function in $\mathbb{C}^m$ such that
\[\limsup_{r \to \infty} \frac{\log T(r,f)}{r} = 0.\]

Assume that $f(z)$ satisfies the differential-difference equation $B(f(z)) Q(f(z))=P(f(z))$,
where $P(f(z))$, $Q(f(z))$ and $B(f(z))$ are defined as in (\ref{cl1}), (\ref{cl2}) and (\ref{cl3}) respectively.
If $\deg(P(f)) \leq n=\deg(B(f))$, then $\parallel\;m(r, Q(f))=o(T(r, f))$.
\end{lem}

\section{\bf{Proofs of Theorems \ref{t2.1}-\ref{t2.3}}}

\subsection {{\bf Proof of Theorem \ref{t2.1}}} 
\begin{proof} Let $f(z)$ be finite order non-constant entire solution of the equation (\ref{pds1}). Then using Lemma \ref{L.3} to (\ref{pds1}), we get
\bea\label{sm.0} \parallel\;m_1T(r,F)+o(T(r,F))=m_2T(r,f(z+c))+o(T(r,f(z+c))).\eea

Let
\bea\label{sdd.2} h(z)=\frac{F^{m_1}(z)-1}{F^{m_1}(z)}.\eea

Clearly $h(z)$ is a non-constant meromorphic function in $\mathbb{C}^m$. Using Lemma \ref{L.3} to (\ref{sdd.2}), we get 
\bea\label{sm.00}\parallel T(r,h)+o(T(r,h))=m_1T(r,F)+o(T(r,F)).\eea

Now from (\ref{pds1}) and (\ref{sdd.2}), it is easy to deduce that 
\[\parallel\;\ol N(r,h)\leq \ol N(r,0, F^{m_1})=\ol N(r,0,F)+o(T(r,F)),\]
\[\parallel\;\ol N(r,0,h)=\ol N(r,1,F^{m_1})\leq \ol N(r,0,f^{m_2}(z+c))=\ol N(r,0,f(z+c))+o(T(r,f(z+c)))\]
and $\parallel \ol N(r,1,h)=0$. Therefore in view of first main theorem and using Lemma \ref{L.2}, we get
\beas \parallel\;m_1T(r,F)&=&T(r,h)+o(T(r,h))\\&\leq& \ol N(r,h)+\ol N(r,0,h)+\ol N(r,1,h)+o(T(r,h))\nonumber\\&\leq&
\ol N(r,0,F)+\ol N(r,0,f(z+c))+o(T(r,F))+o(T(r,f(z+c)))\\&\leq&
T(r,F)+T(r,f(z+c))+o(T(r,F))+o(T(r,f(z+c)))
\eeas
and so from (\ref{sm.0}), we get $\parallel\;\left(m_1-1-m_1/m_2\right)T(r,F)\leq o(T(r,F))$,
which shows that
\bea\label{sm.1}\label{smm.4} 1/m_1+1/m_2\geq 1.\eea

We consider following two cases.\par

\smallskip
{\bf Case 1.} Let $m_1=m_2$. Then (\ref{sm.1}) gives $m_1=m_2=2$ and so $F^2(z)+f^2(z+c)=1$,
i.e.,
\bea\label{sm.2} (F(z)+\iota f(z+c))(F(z)-\iota f(z+c))=1.\eea

Clearly from (\ref{sm.2}), we conclude that both the entire functions $F(z)+\iota f(z+c)$ and $F(z)-\iota f(z+c)$ are of finite order having no zeros. Therefore we assume that
\bea\label{sm.3} F(z)+\iota f(z+c)=e^{\iota \tilde P(z)}\eea
and 
\bea\label{sm.4} F(z)-\iota f(z+c)=e^{-\iota \tilde P(z)},\eea
where $\tilde P(z)$ is a polynomial in $\mathbb{C}^m$. Solving (\ref{sm.3}) and (\ref{sm.4}), we get
\bea\label{sm.5} F(z)=\sideset{}{_{j=1}^m}{\sum} \frac{\partial f(z_1, z_2, \ldots, z_m)}{\partial z_j}=\frac{e^{\iota \tilde P(z)}+e^{-\iota \tilde P(z)}}{2}\eea
and 
\bea\label{sm.6} f(z+c)=\frac{e^{\iota \tilde P(z)}-e^{-\iota \tilde P(z)}}{2\iota}=\sin \tilde P(z).\eea

Obviously $\tilde P(z)$ is a non-constant polynomial. Now from (\ref{sm.6}), we have
\bea\label{sm.7} \frac{\partial f(z_1+c_1,z_2+c_2,\ldots,z_m+c_m)}{\partial z_j}=\frac{e^{\iota \tilde P(z)}+e^{-\iota \tilde P(z)}}{2}\frac{\partial \tilde P(z)}{\partial z_j},\eea
for $j=1,\ldots, m$. Clearly from (\ref{sm.5}) and (\ref{sm.7}), we have
\beas \frac{e^{\iota \tilde P(z+c)}+e^{-\iota \tilde P(z+c)}}{2}=\frac{e^{\iota \tilde P(z)}+e^{-\iota \tilde P(z)}}{2}\sideset{}{_{j=1}^m}{\sum} \frac{\partial \tilde P(z)}{\partial z_j},\eeas
i.e.,
\bea\label{sm.8} \sideset{}{_{j=1}^m}{\sum} \frac{\partial \tilde P(z)}{\partial z_j}e^{\iota (\tilde P(z)+\tilde P(z+c))}+\sideset{}{_{j=1}^m}{\sum}\frac{\partial \tilde P(z)}{\partial z_j}e^{-\iota (\tilde P(z)-\tilde P(z+c))}-e^{2\iota \tilde P(z+c)}=1,\eea
from which we conclude that $\sum_{j=1}^m \frac{\partial \tilde P(z)}{\partial z_j}\not\equiv 0$.  Since $\tilde P(z)$ is a non-constant polynomial, $\tilde P(z)+\tilde P(z+c)$ is also a non-constant polynomial. Now using Lemma \ref{L.6} to (\ref{sm.8}), we get
\bea\label{sm.9} \sideset{}{_{j=1}^m}{\sum} \frac{\partial \tilde P(z)}{\partial z_j}=e^{\iota (\tilde P(z)-\tilde P(z+c))}.\eea

If $\tilde P(z)-\tilde P(z+c)$ is non-constant, then from (\ref{sm.9}) we immediately get a contradiction. Hence $\tilde P(z)-\tilde P(z+c)$ is a constant and so we may assume that 
\bea\label{lm1} \tilde P(z)=d_1z_1+d_2z_2+\ldots+d_mz_m+\hat P(z),\eea
where $d_1,d_2,\ldots,d_m\in \mathbb{C}$ and $\hat P(z)$ is a polynomial in $\mathbb{C}^m$ such that $\hat P(z+c)\equiv \hat P(z)$.

Consequently from (\ref{sm.9}) and (\ref{lm1}), we get
\bea\label{lm2}e^{\iota (d_1c_1+d_2c_2+\ldots+d_mc_m)}\sideset{}{_{j=1}^m}{\sum} \left(\frac{\partial \hat P(z)}{\partial z_j}+d_j\right)=1.\eea
Finally from (\ref{sm.6}), (\ref{lm1}) and (\ref{lm2}), we may assume that 
\[f(z)=\sin (d_1z_1+d_2z_2+\ldots+d_mz_m+P(z)),\]
where $d_1,d_2,\ldots,d_m\in \mathbb{C}$ and $P(z)$ is a polynomial in $\mathbb{C}^m$ such that $P(z+c)\equiv P(z)$ and
\beas e^{\iota (d_1c_1+d_2c_2+\ldots+d_mc_m)}\sideset{}{_{j=1}^m}{\sum} \left(\frac{\partial P(z)}{\partial z_j}+d_j\right)=1.\eeas

\smallskip
{\bf Case 2.} Let $m_1\neq m_2$. We now discuss two sub-cases: one, when $m_2>m_1$ and other, when  $m_1>m_2$.
We now consider following two sub-cases.

\smallskip
{\bf Sub-case 2.1.} Let $m_2>m_1$. Since $\rho(f)<+\infty$, we have 
\[\underset{r \rightarrow \infty}{\limsup } \frac{\log T(r,f)}{r}=0.\]

Then by Lemma \ref{L.7}, we have 
\bea\label{sm.11}\parallel\; T(r,f(z+c))=T(r,f)+o(T(r,f)).\eea

Now using Lemmas \ref{L.1}, \ref{L.3} and (\ref{sm.11}) to (\ref{pds1}), we get
\bea\label{ml0}\parallel\;m_2 T(r,f))=T\left(r, f^{m_2}(z+c)\right)+o(T(r,f))
&=&T\left(r,F^{m_1}-1\right)+o(T(r,f))\\
&=&m_1 T(r,F)+o(T(r,f))\nonumber\\
&=&m_1 m(r,F)+o(T(r,f))\nonumber\\
&\leq & m_1\left(m\left(r,\frac{F}{f}\right)+m(r,f)\right)+o(T(r,f))\nonumber \\
& =&m_1 T(r,f)+o(T(r,f))\nonumber\eea
and so $\parallel (m_2-m_1) T(r,f) \leq o(T(r,f))$, which is impossible.\par

\smallskip
{\bf Sub-case 2.2.} Let $m_1>m_2$. If $m_2\geq 2$, then from (\ref{sm.1}), we get a contradiction. Hence $m_2=1$ and so from (\ref{pds1}), we get
\bea\label{sm.12} F^{m_1}(z)+ f(z+c)= 1.\eea

Differentiating (\ref{sm.12}) partially with respect to $z_1,z_2,\ldots,z_m$ respectively and then adding, we get
\bea\label{sm.13} m_1 F^{m_1-1} \sideset{}{_{j=1}^m}{\sum}\frac{\partial F(z)}{\partial z_j} + \sideset{}{_{j=1}^m}{\sum}\frac{\partial f(z + c)}{\partial z_j} = 0.\eea

Since $F(z)=\sum_{j=1}^m\frac{\partial f(z)}{\partial z_j}$, from (\ref{sm.13}), we get
\bea\label{sm.14} m_1 F^{m_1-1}(z) \sideset{}{_{j=1}^m}{\sum}\frac{\partial F(z)}{\partial z_j} + F(z+c)=0.\eea

We know that $m_1\geq 2$. Now applying Lemma \ref{L.4} to (\ref{sm.14}), we get
\[\parallel m \left(r, \sideset{}{_{j=1}^m}{\sum}\frac{\partial F(z)}{\partial z_j}\right)=o(T(r, F)).\]

Since $F$ is an entire function, we have 
\[\parallel T \left(r, \sideset{}{_{j=1}^m}{\sum}\frac{\partial F(z)}{\partial z_j} \right)=o(T(r, F)).\]

Let $\sum_{j=1}^m\frac{\partial F(z)}{\partial z_j}=\alpha(z)$. Clearly $\alpha(z)$ is an entire function in $\mathbb{C}^m$. Now from (\ref{sm.14}), we get
\bea\label{sm.15} m_1 F^{m_1-1}(z)\alpha(z)=-F(z + c).\eea

If possible suppose that $\alpha\equiv 0$. Then from (\ref{sm.15}), we have $F(z+c)\equiv 0$ and so from (\ref{sm.12}), we get $f(z+2c)\equiv 1$, which is impossible. Hence $\alpha\not\equiv 0$.

Now we consider following two sub-cases.\par

\smallskip
{\bf Sub-case 2.2.1.} Let $m_1>2$. Now applying Lemma \ref{L.4} to (\ref{sm.15}), we get $\parallel m(r,\alpha F)=o(T(r,F))$.
Since $F$ is an entire function and $\parallel T(r,\alpha)=o(T(r,F))$, we have 
\[ T(r,F)=m(r,F)\leq m(r,\alpha F)+m\left(r,0;\alpha\right)\leq o(T(r,F))\]
and so from (\ref{sm.0}), we have $T(r,f(z+c))=o(T(r,f(z+c))$, which contradicts (\ref{sm.11}).\par

\smallskip
{\bf Sub-case 2.2.2.} Let $m_1=2$. Now from equation (\ref{sm.15}), we get
\bea\label{sm.16} 2F(z)\alpha(z)=-F(z+c).\eea

Differentiating (\ref{sm.16}) partially with respect to $z_1,z_2,\ldots,z_m$ respectively and then adding, we get
\bea\label{sm.18} 2\sideset{}{_{j=1}^m}{\sum}\frac{\partial F(z)}{\partial z_j}\alpha(z)+2F(z)\sideset{}{_{j=1}^m}{\sum}\frac{\partial \alpha(z)}{\partial z_j}=-\sideset{}{_{j=1}^m}{\sum}\frac{\partial F(z + c)}{\partial z_j}.\eea

Since $\sum_{j=1}^m\frac{\partial F(z)}{\partial z_j}=\alpha(z)$, from (\ref{sm.18}), we get
\bea\label{sm.19} 2\alpha^2(z)+2F(z)\sideset{}{_{j=1}^m}{\sum}\frac{\partial \alpha(z)}{\partial z_j}=-\alpha(z + c).\eea

If $\sum_{j=1}^m\frac{\partial \alpha(z)}{\partial z_j} \not\equiv 0$, then applying Lemma \ref{L.3} to (\ref{sm.19}), we get $\parallel T(r, F)=o(T(r, F))$ and so from (\ref{sm.0}), we get $T(r,f(z+c))=o(T(r,f(z+c))$, which contradicts (\ref{sm.11}).
Hence $\sum_{j=1}^m\frac{\partial \alpha(z)}{\partial z_j} \equiv 0$ and so from (\ref{sm.19}), we have
\bea\label{sm.20} 2\alpha^2(z)=-\alpha(z+c).\eea

If $\alpha(z)$ is transcendental, then applying Lemma \ref{L.4} to (\ref{sm.20}), we get $\parallel m(r,\alpha)=o(T(r, \alpha))$. This implies $\parallel T(r, \alpha)=o((r, \alpha))$, leading to a contradiction. Hence $\alpha(z)$ is a polynomial in $\mathbb{C}^m$.
It is easy to deduce from (\ref{sm.20}) that $\alpha(z)$ is a constant, say $k(\neq 0)$. Then from (\ref{sm.20}), we have 
$2k^2=-k$, which implies that $k=-\frac{1}{2}$. Therefore, we have 
\bea\label{sm.21} \alpha(z)=\sideset{}{_{j=1}^m}{\sum}\frac{\partial F(z)}{\partial z_j}=-\frac{1}{2}.\eea 

Again from (\ref{sm.16}), we get $F(z)=F(z+c)$, which shows that $F(z)$ is a periodic function with period $c$. 

The Lagrange's auxiliary equations corresponding to (\ref{sm.21}) are as follows
\bea\label{sm.21a} \frac{d z_1}{1}=\frac{dz_2}{1}=\ldots=\frac{dz_m}{1}=\frac{dF(z)}{-\frac{1}{2}}=\frac{dz_1+dz_2+\ldots+dz_m}{m}.\eea

Clearly $z_2-z_1=d_1$, $z_3-z_2=d_2$, $\ldots$, $z_m-z_{m-1}=d_{m-1}$ and $mF(z)+\frac{1}{2}(z_1+z_2+\ldots+z_m)=d_m$ be independent solutions of (\ref{sm.21a}). Therefore the general solution of (\ref{sm.21}) is as follows
\bea\label{sm.22} F(z)=-(z_1+z_2+\ldots+z_m)/2m+g_1(z),\eea
where $g_1(z)=g_1(z_2-z_1,z_3-z_2,\ldots,z_m-z_{m-1})$ is an arbitrary finite order transcendental entire function.
Since $F(z+c)=F(z)$, from (\ref{sm.22}), we have
\[-(z_1 + c_1+z_2+c_2+\ldots+z_m+c_m)/2m + g_1(z+c)=-(z_1+z_2+\ldots+z_m)/2m + g_1(z),\]
which implies
\bea\label{sm.23} g_1(z+c)-g_1(z)=(c_1+c_2+\ldots+c_m)/2m.\eea

Since $F(z+c)=F(z)$, from (\ref{sm.12}), we get $f(z+2c)=1-{F(z+c)}^2=1-F^2(z)=f(z+c)$, which shows that $f(z)$ is $c$-periodic. Then from (\ref{pds1}) and (\ref {sm.22}), we get
\[f(z)=1-\left(-(z_1+z_2+\ldots+z_m)/2m+g_1(z)\right)^2,\]
where $g_1(z)=g_1(z_2-z_1,z_3-z_2,\ldots,z_m-z_{m-1})$ is an arbitrary finite order transcendental entire function such that (\ref{sm.23}) holds. Hence the proof is complete.
\end{proof}

\subsection {{\bf Proof of Theorem \ref{t2.2}}} 
\begin{proof} If possible suppose $f(z)$ is the finite order non-constant meromorphic solution of the equation (\ref{pds1}).
Now from (\ref{lm0}), it is easy to verify that
\bea\label{ml1} N(r,F)\leq N(r,f)+\ol N(r,f)\leq 2 N(r,f).\eea

We divide following two cases.\par

\smallskip
{\bf Case 1.} Let $m_2>m_1\geq 3$. Now from (\ref{pds1}) and (\ref{sdd.2}), it is easy to deduce that 
\[\parallel\;\ol N(r,h)\leq \ol N(r,0, F^{m_1})=\ol N(r,0,F)+o(T(r,F)),\]
\[\parallel\;\ol N(r,0,h)=\ol N(r,1,F^{m_1})\leq \ol N(r,0,f^{m_2}(z+c))=\ol N(r,0,f(z+c))+o(T(r,f(z+c)))\]
and 
\[\parallel\;\ol N(r,1,h)\leq \ol N(r,F^{m_1})=\ol N(r,F)+o(T(r,F)).\] 

Therefore in view of first main theorem and using Lemma \ref{L.2}, we get
\bea\label{sdd.3} \parallel\;m_1T(r,F)&=&T(r,h)+o(T(r,h))\\&\leq& \ol N(r,h)+\ol N(r,0,h)+\ol N(r,1,h)+o(T(r,h))\nonumber\\&\leq&
\ol N(r,0,F)+\ol N(r,0,f(z+c))+\ol N(r,F)+o(T(r,F))+o(T(r,f(z+c)))\nonumber\\&\leq&
T(r,F)+T(r,f(z+c))+o(T(r,F))+o(T(r,f(z+c))).\nonumber
\eea

Now using (\ref{sm.00}) to (\ref{sdd.3}), we get
\beas \parallel\;\left(m_1-2-m_1/m_2\right)T(r,F)\leq o(T(r,F))\eeas
from which we conclude that $m_1-2-m_1/m_2\leq 0$, i.e., $1-2/m_1-1/m_2\leq 0$, i.e., 
\bea\label{smm.27}2/m_1+1/m_2\geq 1.\eea

Since $m_2>m_1\geq 3$, from (\ref{smm.27}), we get a contradiction.\par

\smallskip
{\bf Case 2.} Let $m_2>2m_1$. Now using Lemma \ref{L.1}, we conclude from (\ref{ml0}) and (\ref{ml1}) that
\beas\label{d3}\parallel\;m_2 T(r,f))&=&m_1 T(r,F)+o(T(r,f)) \nonumber\\
& =&m_1 \left(m(r,F)+N(r,F)\right)+o(T(r,f)) \nonumber\\
& \leq & m_1\left(m\left(r,\frac{F}{f}\right)+m(r,f)\right)+2m_1N(r,f)+o(T(r,f))\nonumber \\
& =&2m_1 T(r,f)+o(T(r,f))\eeas
and so $\parallel (m_2-2m_1) T(r,f) \leq o(T(r,f))$, which is impossible.
Hence the proof is complete.
\end{proof}

\subsection {{\bf Proof of Theorem \ref{t2.3}}} 
\begin{proof} Let $f(z)$ be finite order non-constant meromorphic solution of the equation (\ref{pds1}). By the given condition, we see that $F(z)$ and $\frac{\partial f(z)}{\partial z_1}$ share $\infty$ CM. Clearly from (\ref{pds1}), we conclude that $F(z)$ and $f(z+c)$ share $\infty$ IM and so
\bea\label{smm.0}\ol N(r,F)=\ol N(r,f)=\ol N(r,f(z+c)).\eea

Since $\rho(f)<+\infty$, we have $\underset{r \rightarrow \infty}{\limsup } \frac{\log N(r,f)}{r}=0$.
Therefore, we have 
\bea\label{smm.1} \parallel\;N(r,f(z+c))=N(r,f)+o(T(r,f))\eea
(see \cite[(26), pp. 779]{CX1}). Since $F(z)$ and $\frac{\partial f(z)}{\partial z_1}$ share $\infty$ CM, from (\ref{pds1}), we deduce that
\[m_1\left(N(r,f)+\ol N(r,f)\right)=m_2 N(r,f(z+c))\]
and so from (\ref{smm.1}), we get
\bea\label{smm.2} \parallel\;m_1\left(N(r,f)+\ol N(r,f)\right)=m_2 N(r,f)+o(T(r,f)).\eea

Since $m_1\geq m_2$, from (\ref{smm.2}), we get
\bea\label{smm.3} \parallel\;\ol N(r,f)=o(T(r,f)).\eea

If $h(z)=\frac{F^{m_1}(z)-1}{F^{m_1}(z)}$, then using (\ref{pds1}), (\ref{smm.0}) and (\ref{smm.3}), we get
\[\parallel\;\ol N(r,1,h)\leq \ol N(r,F^{m_1})=\ol N(r,F)+o(T(r,F))=o(T(r,f)).\] 

Now proceeding in the same way as done in the proof of Theorem \ref{t2.1}, we arrive at (\ref{sm.1}).

We consider following two cases.\par

\smallskip
{\bf Case 1.} Let $m_1=m_2$. Then (\ref{smm.4}) gives $m_1=m_2=2$ and so $F^2(z)+f^2(z+c)=1$, i.e., $(F(z)+\iota f(z+c))(F(z)-\iota f(z+c))=1$ and so we may assume that
\bea\label{smm.6}F(z)+\iota f(z+c)=h(z)\eea
and
\bea\label{smm.7} F(z)-\iota f(z+c)=1/h(z),\eea
where $h(z)$ is a non-zero meromorphic function on $\mathbb{C}^m$. Solving (\ref{smm.6}) and (\ref{smm.7}), we get 
\bea\label{smm.8} f(z+c)=\frac{h(z)-1/h(z)}{2 \iota}\eea
and 
\bea\label{smm.9} F(z)=\sideset{}{_{j=1}^m}{\sum}\frac{\partial f(z)}{\partial z_j}=\frac{h(z)+1/h(z)}{2}.\eea

Clearly from (\ref{smm.8}) and (\ref{smm.9}), we can easily deduce that
\beas \frac{h(z+c)+1/h(z+c)}{2}=\frac{1}{2 \iota}\left(\sideset{}{_{j=1}^m}{\sum}\frac{\partial h(z)}{\partial z_j}+\frac{\sideset{}{_{j=1}^m}{\sum}\frac{\partial h(z)}{\partial z_j}}{h^2(z)}\right),\eeas
i.e.,
\bea\label{smm.10} \iota\left(h(z+c)+1/h(z+c)\right)=\sideset{}{_{j=1}^m}{\sum}\frac{\partial h(z)}{\partial z_j}\left(1+1/h^2(z)\right).\eea

Applying Lemma \ref{L.3} to (\ref{smm.8}), we get $\parallel T(r,f(z+c))=2T(r,h(z))+o(T(r,h(z))$ and so from (\ref{sm.11}), we have 
\bea\label{smm.11} \parallel\;T(r,f)+o(T(r,f))=2T(r,h)+o(T(r,h)).\eea

Now from (\ref{smm.8}), we obtain
\[\parallel\;\ol N(r,h(z-c))+\ol N(r,0;h(z-c))\leq \ol N(r,f(z))\]
and so from (\ref{smm.3}) and (\ref{smm.11}), we get
\bea\label{smm.10a}\parallel\;\ol N(r,h(z-c))+\ol N(r,0;h(z-c))=o(T(r,h)).\eea

Finally from (\ref{smm.8}), we have $f(z)=\frac{h(z-c)-1/h(z-c)}{2\iota}$, where $h(z)$ is a non-constant finite order meromorphic function in $\mathbb{C}^m$ such that (\ref{smm.10}) and (\ref{smm.10a}) hold.

In particular we consider $(i))$ $\parallel\;N_{(2}(r,f)=o(T(r,f))$ and $(ii))$ $c=0$.

First we suppose that $\parallel\;N_{(2}(r,f)=o(T(r,f))$. Then from (\ref{smm.8}) we obtain
\[\parallel\;N(r,h(z))+N\left(r,0;h(z)\right)\leq N(r,f(z+c))\]
and so from (\ref{smm.1}), we get
\bea\label{smm.12} \parallel\;N(r,h(z))+N(r,0;h(z))\leq N(r,f(z+c))=N(r,f)+o(T(r,f)).\eea
 
Since $\parallel N_{(2}(r,f)=o(T(r,f))$, from (\ref{smm.3}) and (\ref{smm.12}), we obtain
\beas \parallel\;N(r,h(z))+N(r,0;h(z))\leq o(T(r,f))\eeas
and so (\ref{smm.11}) gives
\bea\label{smm.13} \parallel\; N(r,h)+N(r,0;h)=o(T(r,h)).\eea

On the other hand from (\ref{smm.10}), we have 
\bea\label{smm.14} h^2(z)\beta(z)=h(z+c)\sideset{}{_{j=1}^m}{\sum}\frac{\partial h(z)}{\partial z_j},\eea
where
\bea\label{smm.15} \beta(z)=\iota h^2(z+c)-h(z+c)\sideset{}{_{j=1}^m}{\sum}\frac{\partial h(z)}{\partial z_j}+\iota.\eea

Now using Lemma \ref{L.4} to (\ref{smm.14}), we get
$\parallel m(r,\beta)=o(T(r,h)$ and so (\ref{smm.13}) gives $\parallel T(r,\beta)=o(T(r,h))$.
If possible suppose that $\beta\not\equiv \iota.$ Then from (\ref{smm.15}), we get
\[\frac{\beta(z)-\iota}{h^2(z)}=\iota \left(\frac{h(z+c)}{h(z)}\right)^2-\frac{h(z+c)}{h(z)}\frac{\sideset{}{_{j=1}^m}{\sum}\frac{\partial h(z)}{\partial z_j}}{h(z)}\]
and so by Lemmas \ref{L.1} and \ref{L.8}, we get $\parallel m\left(r,\frac{1}{h^2}\right)=o(T(r,h))$. Therefore from (\ref{smm.13}), we get $\parallel T\left(r,\frac{1}{h^2}\right)=o(T(r,h))$ and so by using first main theorem, we obtain $\parallel T(r,h)=o(T(r,h))$, which is impossible. Hence $\beta\equiv \iota$. Clearly from (\ref{smm.15}), we get 
\[\iota h(z+c)=\sideset{}{_{j=1}^m}{\sum}\frac{\partial h(z)}{\partial z_j}.\]

Consequently from (\ref{smm.14}), we have 
\bea\label{smm.16} h^2(z)=h^2(z+c)=-\left(\sideset{}{_{j=1}^m}{\sum}\frac{\partial h(z)}{\partial z_j}\right)^2.\eea

Again from (\ref{smm.8}) and (\ref{smm.9}), we have
\bea\label{smmm.1} \frac{\partial f(z+c)}{\partial z_1}=\frac{\frac{\partial h(z)}{\partial z_1}+\frac{\partial h(z)}{\partial z_1}/h^2(z)}{2 \iota}\eea
and 
\bea\label{smmm.2} F(z+c)=\sideset{}{_{j=1}^m}{\sum}\frac{\partial f(z+c)}{\partial z_j}=\frac{h(z+c)+1/h(z+c)}{2}.\eea

Since $F(z)$ and $\frac{\partial f(z)}{\partial z_1}$ share $\infty$ CM, it follows that $F(z+c)$ and $\frac{\partial f(z+c)}{\partial z_1}$ share $\infty$ CM. On the other hand from (\ref{smm.16}), we see that $h(z)$ and $h(z+c)$ share $0$ and $\infty$ CM. Take $z_0\in\mathbb{C}^m$. Let $\mu^{\infty}_{h(z)}(z_0)>0$. Clearly from (\ref{smmm.2}), we see that $\mu^{\infty}_{F(z+c)}(z_0)=\mu^{\infty}_{h(z+c)}(z_0)=\mu^{\infty}_{h(z)}(z_0)$. Since $F(z+c)$ and $\frac{\partial f(z+c)}{\partial z_1}$ share $\infty$ CM, it follows that 
\[\mu^{\infty}_{\frac{\partial f(z+c)}{\partial z_1}}(z_0)=\mu^{\infty}_{F(z+c)}(z_0)=\mu^{\infty}_{h(z)}(z_0).\]

But from (\ref{smmm.1}), we obtain
\[\mu^{\infty}_{\frac{\partial f(z+c)}{\partial z_1}}(z_0)=\mu^{\infty}_{h(z)}(z_0)+1.\]

Therefore we get a contradiction. Hence $\mu^{\infty}_{h(z)}(z_0)=0$ and so $h(z)$ is an entire function. Similarly if we take $\mu^{0}_{h(z)}(z_0)>0$, then proceeding in the same way, we can prove that $\mu^{0}_{h(z)}(z_0)=0$ and so $h(z)$ has no zeros. Consequently $h$ is an entire function having no zeros. Let $h(z)=e^{\iota \gamma(z)}$, where $\gamma(z)$ is a non-constant polynomial in $\mathbb{C}^m$. Then from (\ref{smm.16}), we get 
\bea\label{d.p1}\sideset{}{_{j=1}^m}{\sum}\frac{\partial \gamma(z)}{\partial z_j}=\pm e^{\iota (\gamma(z+c)-\gamma (z))}.\eea

If $\gamma(z)-\gamma(z+c)$ is non-constant, then from (\ref{d.p1}), we immediately get a contradiction. Hence $\gamma(z)-\gamma(z+c)$ is a constant and so we may assume that 
\bea\label{d.p2} \gamma(z)=d_1z_1+d_2z_2+\ldots+d_mz_m+ P(z),\eea
where $d_1,d_2,\ldots,d_m\in \mathbb{C}$ and $P(z)$ is a polynomial in $\mathbb{C}^m$ such that $P(z+c)\equiv  P(z)$.

Consequently from (\ref{d.p1}) and (\ref{d.p2}), we get
\bea\label{d.p3}e^{-\iota (d_1c_1+d_2c_2+\ldots+d_mc_m)}\sideset{}{_{j=1}^m}{\sum} \left(\frac{\partial P(z)}{\partial z_j}+d_j\right)=\pm 1.\eea
 
Finally from (\ref{smm.8}), (\ref{d.p1}) and (\ref{d.p3}), we may assume that 
\[f(z)=\sin (d_1z_1+d_2z_2+\ldots+d_mz_m+P(z)),\]
where $d_1,d_2,\ldots,d_m\in \mathbb{C}$ and $P(z)$ is a polynomial in $\mathbb{C}^m$ such that $P(z+c)\equiv P(z)$ and
\beas e^{-\iota (d_1c_1+d_2c_2+\ldots+d_mc_m)}\sideset{}{_{j=1}^m}{\sum} \left(\frac{\partial P(z)}{\partial z_j}+d_j\right)=\pm 1.\eeas

Next we suppose that $c=0$. Then from (\ref{smm.10}), we have
\bea\label{smm.17} \sideset{}{_{j=1}^m}{\sum}\frac{\partial h(z)}{\partial z_j}=\iota h(z).\eea

In this case also proceeding in the same way as done above, we can prove that $h(z)$ is an entire function having no zeros. If $h(z)=e^{\iota \gamma(z)}$, where $\gamma(z)$ is a non-constant polynomial in $\mathbb{C}^m$, then from (\ref{smm.17}), we get 
\bea\label{d.p4}\sideset{}{_{j=1}^m}{\sum}\frac{\partial \gamma(z)}{\partial z_j}=1.\eea

The Lagrange's auxiliary equations corresponding to (\ref{d.p4}) are as follows
\bea\label{d.p4a} \frac{d z_1}{1}=\frac{dz_2}{1}=\ldots=\frac{dz_m}{1}=\frac{d\gamma (z)}{1}=\frac{dz_1+dz_2+\ldots+dz_m}{m}.\eea

Clearly $z_2-z_1=d_1$, $z_3-z_2=d_2$, $\ldots$, $z_m-z_{m-1}=d_{m-1}$ and $m\gamma(z)-(z_1+z_2+\ldots+z_m)=d_m$ be independent solutions of (\ref{d.p4a}). Therefore the general solution of (\ref{d.p4}) is as follows
\bea\label{d.p4b} \gamma(z)=(z_1+z_2+\ldots+z_m)/m+P(z),\eea
where $P(z)=P(z_2-z_1,z_3-z_2,\ldots,z_m-z_{m-1})$ is an arbitrary polynomial. Consequently from (\ref{d.p4}), we get
$\sum_{j=1}^m\frac{\partial P(z)}{\partial z_j}=0$. 
Finally from (\ref{smm.8}), we can take 
\[f(z)=\sin \left((z_1+z_2+\ldots +z_m)/m+P(z)\right),\]
where $P(z)=P(z_3-z_2,\ldots,z_m-z_{m-1})$ is a polynomial such that $\sum_{j=1}^m\frac{\partial P(z)}{\partial z_j}=0$.\par

\smallskip
{\bf Case 2.} Let $m_1>m_2$. If $m_2\geq 2$, then (\ref{smm.4}) gives a contradiction. Hence $m_2=1$ and so from (\ref{pds1}), we get
\bea\label{smm.18} F^{m_1}(z)+ f(z + c)= 1.\eea

Differentiating (\ref{smm.18}) partially with respect to $z_1,z_2,\ldots,z_m$ respectively and then adding, we get
\bea\label{smm.19} m_1 F^{m_1-1}(z)\sideset{}{_{j=1}^m}{\sum}\frac{\partial F(z)}{\partial z_j} + F(z + c)=0.\eea

We know that $m_1\geq 2$. Now applying Lemma \ref{L.4} to (\ref{smm.19}), we get
\bea\label{smm.20} \parallel\;m \left(r, \sideset{}{_{j=1}^m}{\sum}\frac{\partial F(z)}{\partial z_j}\right)=o(T(r, F).\eea

Clearly from (\ref{smm.2}), we have
\bea\label{smm.21} \parallel\;N(r,f)=o(T(r,f)).\eea

On the other hand from (\ref{sm.0}) and (\ref{sm.11}), we get
\bea\label{smm.22} \parallel\;m_1T(r,F)+o(T(r,F))=T(r,f)+o(T(r,f)).\eea

Now using (\ref{smm.21}) and (\ref{smm.22}), we obtain
\bea\label{smm.23} \parallel\;N\left(r,\sideset{}{_{j=1}^m}{\sum}\frac{\partial F(z)}{\partial z_j}\right)=o(T(r,F)).\eea

Therefore from (\ref{smm.20}) and (\ref{smm.23}), we have 
\[\parallel T \left(r, \sideset{}{_{j=1}^m}{\sum}\frac{\partial F(z)}{\partial z_j} \right)=o(T(r, F)).\]

Let $\sum_{j=1}^m\frac{\partial F(z)}{\partial z_j}=\alpha(z)\not\equiv 0$. Then from (\ref{smm.19}), we get
\bea\label{smm.24} m_1 F^{m_1-1}(z)\alpha(z)=-F(z + c).\eea

If possible suppose that $m_1>2$. Now using Lemma \ref{L.4} to (\ref{smm.24}), we obtain $\parallel m(r,\alpha F)=o(T(r,F)$. Clearly  from (\ref{smm.21}) and (\ref{smm.22}), we have $\parallel N(r,F)=o(T(r,F))$ and so 
\[\parallel\;T(r,F)=m(r,F)+N(r,F)\leq m(r,\alpha F)+m\left(r,0;\alpha\right)\leq o(T(r,F)).\]

Consequently (\ref{smm.22}) yields $\parallel T(r,f)=o(T(r,f)$, which is impossible. Hence $m_1=2$ and so from (\ref{smm.24}), we get
\bea\label{smm.25} 2F(z)\alpha(z)=-F(z + c).\eea

Differentiating (\ref{smm.25}) partially with respect to $z_1,z_2,\ldots,z_m$ respectively and then adding and using $\sum_{j=1}^m\frac{\partial F(z)}{\partial z_j}=\alpha(z)$, we get
\bea\label{smm.26} 2\alpha^2(z)+2F(z)\sideset{}{_{j=1}^m}{\sum}\frac{\partial \alpha(z)}{\partial z_j}=-\alpha(z + c).\eea

If $\sum_{j=1}^m\frac{\partial \alpha(z)}{\partial z_j} \not\equiv 0$, then applying Lemma \ref{L.3} to (\ref{smm.26}), we get $\parallel T(r, F)=o(T(r, F))$, which is impossible. Hence $\sum_{j=1}^m\frac{\partial \alpha(z)}{\partial z_j} \equiv 0$ and so from (\ref{smm.26}), we obtain $2\alpha^2(z)=-\alpha(z + c)$. Now proceeding in the same way as done in the proof of Sub-case 2.2.2 of Theorem \ref{t2.1}, we get the desire conclusion. Hence the proof is complete.
\end{proof}

\subsection {{\bf Proof of Theorem \ref{t2.4}}} 
\begin{proof} Let $f(z)$ be finite order non-constant meromorphic solution of the equation (\ref{pds1}) such that the set of
poles of $f(z)$ is algebraic. Clearly $N(r,f)=O(\log r)$ and $N(r,f)=O(\log r)$.
If $h(z)=\frac{F^{m_1}(z)-1}{F^{m_1}(z)}$, then using (\ref{pds1}), we get
\[\parallel\;\ol N(r,1,h)\leq \ol N(r,F^{m_1})=\ol N(r,F)+o(T(r,F))=o(T(r,f)).\] 

Now proceeding in the same way as done in the proof of Theorem \ref{t2.1}, we arrive at (\ref{sm.1}).

We consider following there cases.\par

\smallskip
{\bf Case 1.} Let $m_1=m_2$. Then (\ref{sm.1}) gives $m_1=m_2=2$ and so $F^2(z)+f^2(z+c)=1$. Since $N_{(2}(r,f)=O(\log r)$, 
proceeding in the same way as done in the proof of Case 1 of Theorem \ref{t2.3}, we can prove that
$f(z)=\sin (d_1z_1+\ldots+d_mz_m+P(z))$,
where $d_1,\ldots,d_m\in \mathbb{C}$ and $P(z)$ is a polynomial in $\mathbb{C}^m$ such that $P(z+c)\equiv P(z)$ and
\beas e^{-\iota (d_1c_1+d_2c_2+\ldots+d_mc_m)}\sideset{}{_{j=1}^m}{\sum} \left(\frac{\partial P(z)}{\partial z_j}+d_j\right)=\pm 1.\eeas

\smallskip
{\bf Case 2.} Let $m_2>m_1$. In this case proceeding in the same way as done in the proof of Sub-case 2.1 of Theorem \ref{t2.1}, we get a contradiction.\par

\smallskip
{\bf Case 2.} Let $m_2>m_1$. Now proceeding in the same way as done in the proof of Case 2 of Theorem \ref{t2.3},  we get the desire conclusion. Hence the proof is complete.
\end{proof}

\vspace{0.1in}
{\bf Compliance of Ethical Standards:}\par

{\bf Conflict of Interest.} The authors declare that there is no conflict of interest regarding the publication of this paper.\par

{\bf Data availability statement.} Data sharing not applicable to this article as no data sets were generated or analysed during the current study.

\end{document}